\documentclass[letterpaper]{amsart}

\usepackage[utf8]{inputenc}
\usepackage{lmodern}
\usepackage[T1]{fontenc}
\usepackage{amssymb}
\usepackage{enumitem}
\usepackage{mathtools}

\usepackage{tikz}
\usetikzlibrary{calc}
\usepackage{tikz-cd}
\usepackage[pdfusetitle]{hyperref}
\usepackage{cleveref}
\usepackage[hyphenbreaks]{breakurl}

\tikzcdset{arrow style=Latin Modern}

\newlist{enumarabic}{enumerate}{1}
\setlist[enumarabic]{font=\normalfont,label=(\arabic*),leftmargin=0.3in}
\newlist{enumroman}{enumerate}{1}
\setlist[enumroman]{font=\normalfont,label=(\roman*),leftmargin=0.3in}

\numberwithin{equation}{section}
\allowdisplaybreaks[4]

\theoremstyle{plain}
\newtheorem{theorem}{Theorem}[section]

\theoremstyle{definition}

\newtheorem{remark}[theorem]{Remark}
\newtheorem{example}[theorem]{Example}

\theoremstyle{remark}

\DeclareMathAlphabet\mathbfit{OML}{cmm}{b}{it}

\let\newterm\emph

\def\arxiv#1{\burlalt{http://arxiv.org/abs/#1}{arXiv:#1}}
\def\doi#1{\burlalt{http://doi.org/#1}{doi:#1}}
\def\doiwithbreak#1#2{\burlalt{http://doi.org/#1}{doi:#2}}

\def\cf{\emph{cf.}}

\let\epsilon\varepsilon
\let\phi\varphi
\let\emptyset\varnothing

\def\bigpair#1{{\bigl\langle#1\bigr\rangle}}

\def\Z{\mathbb{Z}}

\def\R{\mathbb{R}}
\def\C{\mathbb{C}}
\def\CP{\mathbb{CP}}
\def\HT{H_{T}}
\def\HG{H_{G}}
\def\KT{K_{T}}

\def\mm{\mathfrak{m}}
\DeclareMathOperator{\Hom}{Hom}
\DeclareMathOperator{\rank}{rank}
\def\PP{PP}
\def\KK{\mathbf{K}}
\def\TC{T_{\C}}

\def\ww{30}
\def\aa{0mm}
\def\rr{0.075cm}
\def\xr{2.55mm}
\def\yr{1mm}

\begin{document}

\title{The Chang--Skjelbred lemma and generalizations}
\author{Matthias Franz}
\thanks{The author was supported by an NSERC Discovery Grant.}
\address{Department of Mathematics, University of Western Ontario,
  London, Ont.\ N6A\;5B7, Canada}
\email{mfranz@uwo.ca}


\begin{abstract}
  We review the Chang--Skjelbred lemma for torus-equivariant cohomology and discuss several generalizations of it:
  to other coefficients, other groups and also to syzygies in equivariant cohomology and the Atiyah--Bredon sequence.
\end{abstract}

\maketitle

\section{The Chang--Skjelbred lemma}
\label{sec:cs-lemma}

Let \(T\cong(S^{1})^{n}\) be a torus of rank~\(n\ge0\), and let \(X\) be a ``sufficiently nice'' \(T\)-space.
For example, \(X\) can be a compact smooth \(T\)-manifold, possibly with boundary, or a complex algebraic variety with an
action of the complexification~\(\TC=(\C^{\times})^{n}\).

Let \(ET\to BT\) be the universal \(T\)-bundle.
By definition, the \(T\)-equivariant cohomology of~\(X\) is the cohomology of the Borel construction
(or homotopy quotient)~\(X_{T}=ET\times_{T}X\),
\begin{equation}
  \HT^{*}(X) = H^{*}(X_{T}).
\end{equation}
Unless specified otherwise, we take cohomology with coefficients in~\(\R\).

Recall that \(R=H^{*}(BT)\) is a polynomial ring in generators~\(t_{1}\),~\dots,~\(t_{n}\) of degree~\(2\).
Via the canonical projection \(X_{T}\to BT\), \(\HT^{*}(X)\) becomes a module (even an algebra) over~\(R\).
Moreover, under our assumptions on~\(X\), \(\HT^{*}(X)\) is a finitely generated \(R\)-module.
If \(T\) acts trivially on~\(X\), then \(X_{T}=BT\times X\), so that we get an isomorphism of \(R\)-algebras~\(\HT^{*}(X)\cong R\otimes H^{*}(X)\).

The \(T\)-space~\(X\) is said to be \newterm{equivariantly formal} if \(\HT^{*}(X)\) is a free \(R\)-module.
There are several equivalent conditions for equivariant formality:
\begin{enumarabic}
\item Serre spectral sequence for the bundle~\(X\hookrightarrow X_{T}\to BT\)
  degenerates on the second page~\(E_{2}=R\otimes H^{*}(X)\).
\item There is an isomorphism of \(R\)-modules
  \begin{equation}
    \HT^{*}(X) \cong R\otimes H^{*}(X).
  \end{equation}
  (which is \emph{not} multiplicative in general).
\item The canonical restriction map~\(\HT^{*}(X)\to H^{*}(X)\) induced by the inclusion of the fibre is surjective.
  (In traditional terms, ``\(X\) is totally non-homologous to~\(0\) in~\(X_{T}\)''.)
\item There is an isomorphism of graded algebras
  \begin{equation}
    \label{eq:H-from-HT}
    H^{*}(X) \cong \HT^{*}(X) \bigm/ \mm\cdot\HT^{*}(X),
  \end{equation}
  where \(\mm=(t_{1},\dots,t_{n})\lhd R\) is the maximal homogeneous ideal.
\item The sum of the Betti numbers of the fixed point set~\(X^{T}\) agrees with that of~\(X\),
  \begin{equation}
    \label{eq:betti-sum}
    \dim H^{*}(X^{T}) = \dim H^{*}(X),
  \end{equation}
  \cf~\cite[Thm.~3.10.4]{AlldayPuppe:1993}.
\end{enumarabic}
We additionally list the following sufficient criteria: \(X\) is equivariantly formal if
\begin{enumarabic}[resume]
\item \(H^{*}(X)\) vanishes in odd degrees (because all differentials in the Serre spectral sequence must be zero),
\item \(X\) is a compact Kähler manifold and \(X^{T}\ne\emptyset\) (Blanchard~\cite[Sec.~II.1]{Blanchard:1956}),
\item \(X\) is a compact symplectic manifold with a Hamiltonian \(T\)-action
  (Frankel~\cite{Frankel:1959}, Kirwan~\cite[Prop.~5.8]{Kirwan:1984}), or
\item \(X\) is a complete (for instance, projective) smooth complex algebraic variety
  with an algebraic action of the algebraic torus~\(\TC\)
  (Goresky--Kott\-witz--MacPherson~\cite[Thm.~14.1\,(2)]{GoreskyKottwitzMacPherson:1998},
  Weber~\cite{Weber:2003}).
\end{enumarabic}

The \newterm{equivariant \(1\)-skeleton}~\(X_{1}\subset X\) is the union of the fixed points and the orbits of dimension~\(1\).
In other words, it consists of all~\(x\in X\) whose isotropy group~\(T_{x}\subset T\) has rank at least~\(n-1\).
Note that even for smooth actions, \(X_{1}\) is not smooth in general.
For systematic reasons, we also write \(X_{0}=X^{T}\).

\begin{theorem}[Chang--Skjelbred lemma]
  If \(X\) is equivariantly formal, then the \newterm{Chang--Skjelbred sequence}
  \begin{equation*}
    0 \longrightarrow \HT^{*}(X) \stackrel{\iota^{*}}{\longrightarrow} \HT^{*}(X_{0}) \stackrel{\delta}{\longrightarrow} \HT^{*+1}(X_{1},X_{0})
  \end{equation*}
  is exact. Here \(\iota\colon X_{0}\hookrightarrow X\) is the inclusion map,
  and \(\delta\) is the connecting homomorphism~
  in the long exact sequence for the pair~\((X_{1},X_{0})\).
\end{theorem}

In words, this means that \(\HT^{*}(X)\) embeds into~\(\HT^{*}(X_{0})\) via~\(\iota^{*}\)
and that the image equals the kernel of~\(\delta\).
It entails that the equivariant cohomology of an equivariantly formal \(T\)-space~\(X\) can be computed
out of the equivariant \(1\)-skeleton~\(X_{1}\) alone.
This often provides an efficient way to compute \(\HT^{*}(X)\) and, via the isomorphism~\eqref{eq:H-from-HT}, also \(H^{*}(X)\).
Note that we even obtain the product structure in~\(\HT^{*}(X)\) and~\(H^{*}(X)\) this way because the map~\(\iota^{*}\) is multiplicative
and the image of~\(\HT^{*}(X)\) therefore a subring of~\(\HT^{*}(X_{0})\cong R\otimes H^{*}(X_{0})\).
As mentioned above, the product structure of the latter is componentwise and therefore comparatively easy to understand.

Similarly, if \(f\colon X\to Y\) is a map between equivariantly formal \(T\)-spaces, then the commutative diagram
\begin{equation}
  \begin{tikzcd}
    \HT^{*}(Y) \arrow{d}[left]{f^{*}} \arrow[hook]{r}{\iota^{*}} & \HT^{*}(Y_{0}) \arrow{d}{f^{*}} \\
    \HT^{*}(X) \arrow[hook]{r}{\iota^{*}} & \HT^{*}(X_{0})
  \end{tikzcd}
\end{equation}
allows to reconstruct the induced map~\(f^{*}\colon\HT^{*}(Y)\to\HT^{*}(X)\) from the restriction of~\(f\) to the fixed point sets.

\begin{remark}
  Given the long exact sequence of the pair~\((X_{1},X_{0})\), the exactness of the Chang--Skjelbred sequence
  is equivalent to the condition that the images of~\(\HT^{*}(X)\) and~\(\HT^{*}(X_{1})\) in~\(\HT^{*}(X_{0})\) coincide.
  The latter image moreover is the intersection of all images~\(\HT^{*}(X^{K})\to\HT^{*}(X_{0})\), where \(K\subset T\) runs
  through the (finitely many\footnote{%
    Torus representations have only finitely many orbit types, hence so do compact smooth \(T\)-manifolds
    by the slice theorem. An algebraic \(\TC\)-variety can be decomposed into finitely many smooth \(\TC\)-varieties,
    which according to Sumihiro's theorem can be covered by finitely many affine \(\TC\)-varieties. The latter can be embedded into
    representation spaces, which forces the number of orbit types to be finite.}%
  ) codimension-\(1\) subtori that occur as the identity components of some isotropy group in~\(X\).
  This is because one has an isomorphism
  \begin{equation}
    \HT^{*}(X_{1},X_{0}) = \bigoplus_{K} \HT^{*}(X^{K},X_{0})
  \end{equation}
  where \(K\) runs over the same subtori as before. (Note that \(X^{K}=X_{0}\) for a codimension-\(1\) subtorus
  that does not occur as the identity component of an isotropy group in~\(X\).)
\end{remark}

\begin{remark}
  The Chang--Skjelbred lemma appeared in~\cite[Lemma~2.3]{ChangSkjelbred:1974}
  (with a slightly weaker assumption on~\(\HT^{*}(X)\) than freeness).
  Around the same time, Atiyah established a longer exact sequence
  in the context of equivariant \(K\)-theory, see \Cref{sec:other}.
  Several decades later Goresky--Kottwitz--MacPherson rediscovered the Chang--Skjelbred lemma and popularized it.
  In particular, they pointed out its computational power in the special case where \(X_{0}\) consists of finitely many
  points and \(X_{1}\) of finitely many \(2\)-spheres that are rotated via characters~\(T\to S^{1}\) and glued together at their poles.
  This is now called ``GKM theory'' and has numerous applications in symplectic geometry, algebraic geometry and combinatorics.
  For more about the theory and applications of GKM theory, see the articles~\cite{Holm:2008} and~\cite{Tymoczko:2005}
  as well as the forthcoming papers~\cite{Goldin:tufts} and~\cite{Tymoczko:tufts}.
\end{remark}

\section{Examples}

\begin{example}
  We consider the standard action of~\(T=(S^{1})^{2}\) on~\(X=\CP^{2}\)
  given by~\((g_{1},g_{2})\cdot[x_{0}:x_{1}:x_{2}]=[x_{0}:g_{1}x_{1}:g_{2}x_{2}]\).
  This space is equivariantly formal by any of the three necessary conditions listed in \Cref{sec:cs-lemma}.

  The fixed point set~\(X_{0}\) consists of the three points
  \begin{equation}
    x_{0} = [1,0,0], \quad x_{1} = [0,1,0], \quad x_{2} = [0,0,1]
  \end{equation}
  and the equivariant \(1\)-skeleton~\(X_{1}\) is the union of three \(2\)-spheres,
  \begin{equation}
    S(\chi_{1}) = [*,*,0], \quad S(\chi_{2}) = [*,0,*], \quad S(\chi_{2}/\chi_{1}) = [0,*,*].
  \end{equation}
  Here \(S(\chi)\) denotes a \(2\)-sphere on which \(T\) acts by rotations via the character~\(\chi\colon T\to S^{1}\).
  The characters~\(\chi_{1}\) and~\(\chi_{2}\) are the canonical projections.
  The left picture in \Cref{fig:cp2} shows how these three spheres are glued together at their fixed points.
  The quotient~\(X_{1}/T\) is a graph, displayed next to it.

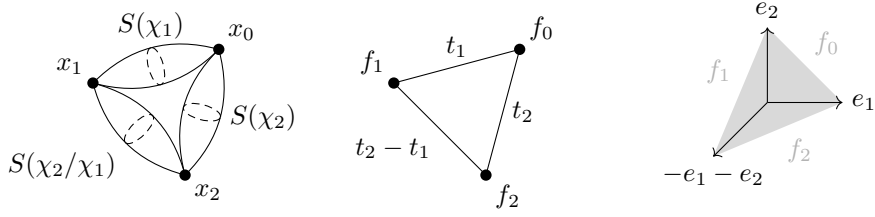
\begin{figure}[ht]
  \begin{tikzpicture}[text opacity=1]
    \coordinate (oo) at (0,0);
    \coordinate (v1) at (165:1);
    \coordinate (v2) at (285:1);
    \coordinate (v3) at (45:1);
    \draw ($(v1)!\aa!-\ww:(v2)$) to [out=-45-\ww, in=135+\ww] node (n12o) {} ($(v2)!\aa!+\ww:(v1)$);
    \draw ($(v1)!\aa!+\ww:(v2)$) to [out=-45+\ww, in=135-\ww] node (n12i) {} ($(v2)!\aa!-\ww:(v1)$);
    \draw ($(v3)!\aa!-\ww:(v1)$) to [out=-165-\ww, in=15+\ww] node (n13o) {} ($(v1)!\aa!+\ww:(v3)$);
    \draw ($(v3)!\aa!+\ww:(v1)$) to [out=-165+\ww, in=15-\ww] node (n13i) {} ($(v1)!\aa!-\ww:(v3)$);
    \draw ($(v3)!\aa!-\ww:(v2)$) to [out=-105-\ww, in=75+\ww] node (n23i) {} ($(v2)!\aa!+\ww:(v3)$);
    \draw ($(v3)!\aa!+\ww:(v2)$) to [out=-105+\ww, in=75-\ww] node (n23o) {} ($(v2)!\aa!-\ww:(v3)$);
    \draw [densely dashed] ($(n23i)!0.5!(n23o)$) circle [x radius=\xr, y radius=\yr, rotate=-15];
    \draw [densely dashed] ($(n12i)!0.5!(n12o)$) circle [x radius=\xr, y radius=\yr, rotate=-135];
    \draw [densely dashed] ($(n13i)!0.5!(n13o)$) circle [x radius=\xr, y radius=\yr, rotate=105];
    \fill (v1) circle [radius=\rr];
    \fill (v2) circle [radius=\rr];
    \fill (v3) circle [radius=\rr];
    \node [above left] at (v1) {$x_{1}$};
    \node [below right] at (v2) {$x_{2}$};
    \node [above right] at (v3) {$x_{0}$};
    \node [below left] at (n12o) {$S(\chi_{2}/\chi_{1})$};
    \node [above] at (n13o) {$S(\chi_{1})$};
    \node [right] at (n23o) {$S(\chi_{2})$};
    \pgftransformxshift{4cm}
    \coordinate (oo) at (0,0);
    \coordinate (v1) at (165:1);
    \coordinate (v2) at (285:1);
    \coordinate (v3) at (45:1);
    \draw (v1) -- node [below left] {$t_{2}-t_{1}$} (v2) -- node [right] {$t_{2}$} (v3) -- node [above] {$t_{1}$} cycle;
    \fill (v1) circle [radius=\rr];
    \fill (v2) circle [radius=\rr];
    \fill (v3) circle [radius=\rr];
    \node [above left] at (v1) {$f_{1}$};
    \node [below right] at (v2) {$f_{2}$};
    \node [above right] at (v3) {$f_{0}$};
    \pgftransformxshift{4cm}
    \coordinate (oo) at (0,0);
    \coordinate (v1) at (1,0);
    \coordinate (v2) at (0,1);
    \coordinate (v3) at (-0.72,-0.72);
    \filldraw [fill=gray, fill opacity=0.3, draw=none] (oo) -- (v1) -- node [above right] {$f_{0}$} (v2);
    \filldraw [fill=gray, fill opacity=0.3, draw=none] (oo) -- (v2) -- node [above left] {$f_{1}$} (v3);
    \filldraw [fill=gray, fill opacity=0.3, draw=none] (oo) -- (v3) -- node [below right] {$f_{2}$} (v1);
    \draw [->] (oo) -- (v1) node [right] {$e_{1}$};
    \draw [->] (oo) -- (v2) node [above] {$e_{2}$};
    \draw [->] (oo) -- (v3) node [below] {$-e_{1}-e_{2}$};
  \end{tikzpicture}
  \caption{Equivariant $1$-skeleton, GKM graph and fan for~$\CP^{2}$}
  \label{fig:cp2}
\end{figure}
  
  Recall that a character~\(\chi\) can be identified with an element of~\(H^{1}(T;\Z)\) and also with a linear polynomial~\(\ell\in H^{2}(BT;\Z)\).
  For~\(S=S(\chi)\) one has under this identification
  \begin{equation*}
    \HT^{*+1}(S,S_{0}) \cong R/(\ell),
  \end{equation*}
  and for a fixed point~\(x\in S_{0}\) the map~\(\delta:\HT^{*}(\{x\})\to\HT^{*+1}(S,S_{0})\)
  corresponds up to sign to the canonical projection~\(R\to R/(\ell)\).
  
  An element of~\(\HT^{*}(X_{0})\) is a triple~\((f_{0},f_{1},f_{2})\in R\oplus R\oplus R\) of polynomials.
  Such a triple is in the kernel of~\(\delta\) if and only if
  \begin{equation}
    \label{eq:cond-mod}
    f_{1}\equiv f_{0} \pmod{t_{2}},
    \qquad
    f_{2}\equiv f_{0} \pmod{t_{1}},
    \qquad
    f_{2}\equiv f_{1} \pmod{t_{2}-t_{1}}.
  \end{equation}
  In the graph above, we have labelled the edges by the polynomials giving the divisibility conditions.
  Such a graph is called a \newterm{GKM graph}; from it one can directly read off the description of~\(\HT^{*}(X)\)
  we have given.

  The projective space~\(X=\CP^{2}\) is actually a complete smooth toric variety, given by the fan~\(\Sigma\)
  on the right of \Cref{fig:cp2}. Hence we see that
  \begin{equation}
    \HT^{*}(X) \cong \PP(\Sigma) = \bigl\{\, f\colon |\Sigma|\to\R \bigm| \text{\(f\) piecewise polynomial}\,\bigr\},
  \end{equation}
  where \(|\Sigma|\) is the support of the fan (equal to~\(\R^{2}\) in this case),
  and ``piecewise polynomial'' means that \(f\) restricts to a polynomial function on each (closed) cone.
  The congruences~\eqref{eq:cond-mod} amount to the condition that the function~\(|\Sigma|\to\R\) given
  by the~\(f_{i}\) 
  is well-defined on the \(1\)-di\-men\-sional pairwise intersections of the \(2\)-dimensional cones.
\end{example}

\begin{example}
  Now consider~\(X=\CP^{2}\setminus\{[1,0,0]\}\). It equivariantly deformation-retracts to~\(\CP^{1}\),
  so that is again equivariantly formal. 
  The equivariant \(1\)-skeleton consists of the sphere~\(\{[0,*,*]\}\) and two punctured spheres.

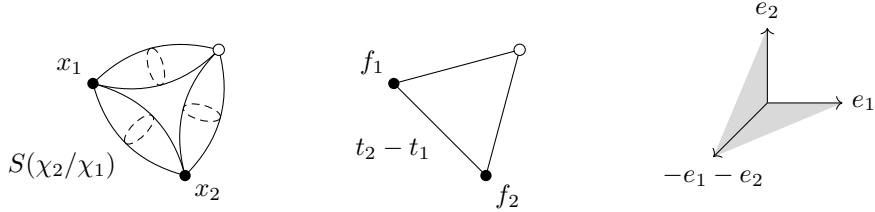
\begin{figure}[h]
  \begin{tikzpicture}
    \coordinate (oo) at (0,0);
    \coordinate (v1) at (165:1);
    \coordinate (v2) at (285:1);
    \coordinate (v3) at (45:1);
    \draw ($(v1)!\aa!-\ww:(v2)$) to [out=-45-\ww, in=135+\ww] node (n12o) {} ($(v2)!\aa!+\ww:(v1)$);
    \draw ($(v1)!\aa!+\ww:(v2)$) to [out=-45+\ww, in=135-\ww] node (n12i) {} ($(v2)!\aa!-\ww:(v1)$);
    \draw ($(v3)!\aa!-\ww:(v1)$) to [out=-165-\ww, in=15+\ww] node (n13o) {} ($(v1)!\aa!+\ww:(v3)$);
    \draw ($(v3)!\aa!+\ww:(v1)$) to [out=-165+\ww, in=15-\ww] node (n13i) {} ($(v1)!\aa!-\ww:(v3)$);
    \draw ($(v3)!\aa!-\ww:(v2)$) to [out=-105-\ww, in=75+\ww] node (n23i) {} ($(v2)!\aa!+\ww:(v3)$);
    \draw ($(v3)!\aa!+\ww:(v2)$) to [out=-105+\ww, in=75-\ww] node (n23o) {} ($(v2)!\aa!-\ww:(v3)$);
    \draw [densely dashed] ($(n23i)!0.5!(n23o)$) circle [x radius=\xr, y radius=\yr, rotate=-15];
    \draw [densely dashed] ($(n12i)!0.5!(n12o)$) circle [x radius=\xr, y radius=\yr, rotate=-135];
    \draw [densely dashed] ($(n13i)!0.5!(n13o)$) circle [x radius=\xr, y radius=\yr, rotate=105];
    \fill (v1) circle [radius=\rr];
    \fill (v2) circle [radius=\rr];
    \filldraw[fill=white] (v3) circle [radius=\rr];
    \node [above left] at (v1) {$x_{1}$};
    \node [below right] at (v2) {$x_{2}$};
    \node [below left] at (n12o) {$S(\chi_{2}/\chi_{1})$};
    \pgftransformxshift{4cm}
    \coordinate (oo) at (0,0);
    \coordinate (v1) at (165:1);
    \coordinate (v2) at (285:1);
    \coordinate (v3) at (45:1);
    \draw (v1) -- node [below left] {$t_{2}-t_{1}$} (v2) -- (v3) -- cycle;
    \fill (v1) circle [radius=\rr];
    \fill (v2) circle [radius=\rr];
    \filldraw[fill=white] (v3) circle [radius=\rr];
    \node [above left] at (v1) {$f_{1}$};
    \node [below right] at (v2) {$f_{2}$};
    \pgftransformxshift{4cm}
    \coordinate (oo) at (0,0);
    \coordinate (v1) at (1,0);
    \coordinate (v2) at (0,1);
    \coordinate (v3) at (-0.72,-0.72);
    \fill [fill=gray, fill opacity=0.3] (oo) -- (v2) -- (v3);
    \fill [fill=gray, fill opacity=0.3] (oo) -- (v3) -- (v1);
    \draw [->] (oo) -- (v1) node [right] {$e_{1}$};
    \draw [->] (oo) -- (v2) node [above] {$e_{2}$};
    \draw [->] (oo) -- (v3) node [below] {$-e_{1}-e_{2}$};
  \end{tikzpicture}
  \caption{Equivariant $1$-skeleton, GKM graph and fan for~$\CP^{2}\setminus\{x_{0}\}$}
\end{figure}

  A punctured sphere can be equivariantly deformation-retracted to the other pole,
  so that the relative equivariant cohomology of the pair vanishes.
  This implies that the two punctured spheres can be ignored
  when computing the kernel of~\(\delta\). Hence \(\HT^{*}(X)\) is isomorphic to the algebra
  of pairs~\((f_{1},f_{2})\in R\oplus R\) satisfying the relation
  \begin{equation}
    f_{2}\equiv f_{1} \pmod{t_{2}-t_{1}}.
  \end{equation}

  Note that we are again looking at a smooth toric variety. Its fan~\(\Sigma\) is obtained from the one for~\(\CP^{2}\)
  by removing the cone covering the first quadrant.
  We therefore find again that \(\HT^{*}(X)\) is isomorphic to the piecewise polynomials on the support of~\(\Sigma\).
  By a result of Brion~\cite[Sec.~2.2]{Brion:1996}, this holds in fact for \emph{all} smooth toric varieties,
  including those which are not equivariantly formal.
\end{example}

\section{Integer coefficients}

Let us indicate what happens if we take cohomology with integer coefficients.
Then \(R=H^{*}(BT)=\Z[t_{1},\dots,t_{n}]\) is still a polynomial ring, but of course with integer coefficients this time.
There are several ways to generalize the notion of equivariant formality to integer coefficients:
\begin{enumarabic}
\item \label{cond1} \(\HT^{*}(X)\) is free over~\(R\).
\item \label{cond2} There is an isomorphism of \(R\)-modules~\(\HT^{*}(X)\cong R\otimes H^{*}(X)\).
\item \label{cond3} The restriction map~\(\HT^{*}(X)\to H^{*}(X)\) is surjective.
\end{enumarabic}
Then \(\hbox{\ref{cond1}}\Rightarrow\hbox{\ref{cond2}}\Rightarrow\hbox{\ref{cond3}}\), but the reverse
implications do not hold in general. (To see that \ref{cond2} does not imply \ref{cond1}, consider
a trivial \(T\)-space with torsion in its cohomology. For an example where \ref{cond3} holds, but not \ref{cond2}
we refer to~\cite[Example~5.2]{FranzPuppe:2007}.) If \(H^{*}(X)\) is free over~\(\Z\), however, then all three
conditions are equivalent.

\begin{theorem}[Franz--Puppe~{\cite[Thm.~1.1]{FranzPuppe:2007}}]
  Assume that condition~\ref{cond3} holds and that all isotropy groups in~\(X\) are connected.
  Then the Chang--Skjelbred sequence is exact over~\(\Z\).
\end{theorem}

The connectivity assumption on the isotropy groups can be weakened, in particular if \ref{cond1} holds,
see~\cite[Cor.~2.2]{FranzPuppe:2011} and also Anderson--Fulton~\cite[Thm.~3.4]{AndersonFulton:2021}.

\section{Reflexive modules}

We again take cohomology with real coefficients, so that \(R=\R[t_{1},\dots,t_{n}]\).

Let \(M\) be a finitely generated \(R\)-module. We write \(M^{\vee}=\Hom_{R}(M,R)\) for the dual module.
Recall that a (graded) commutative bilinear pairing~\(M\times M\to R\) is called \newterm{perfect}
if it identifies \(M\) with~\(M^{\vee}\).

Moreover, \(M\) is called \newterm{reflexive} if the canonical map to its double dual,
\begin{equation}
  \label{eq:def-reflexive}
  M \to M^{\vee\vee},
  \quad
  m \mapsto \bigl( \gamma \mapsto \gamma(m) \bigr),
\end{equation}
is an isomorphism. (This is actually equivalent to \(M\) being isomorphic
to the \(R\)-dual of some finitely generated \(R\)-module~\(N\).) This condition is weaker than freeness,
but stronger than torsion-freeness over~\(R\). For example, the maximal homogeneous ideal~\(\mm\lhd R\)
is torsion-free (like any other submodule of a free \(R\)-module), but it is not reflexive for~\(n\ge2\):
The restriction map~\(R\cong R^{\vee}\to\mm^{\vee}\) is an isomorphism in this case, hence \(\mm^{\vee\vee}\cong R\)
with the map~\eqref{eq:def-reflexive} corresponding to the canonical inclusion~\(\mm\hookrightarrow R\).
For~\(n\ge3\), the second syzgygy of the Koszul resolution~\(\KK^{*}\) of~\(\R\) over~\(R\)
(that is, the kernel of the map~\(\KK^{-2}\to\KK^{-1}\)) is an example of a reflexive \(R\)-module that is not free.

\begin{theorem}[Allday--Franz--Puppe~{\cite[Thm.~1.1, Cor.~1.3]{AlldayFranzPuppe:orbits1}}]
  \label{thm:reflexive}
  The following two conditions are equivalent:
  \begin{enumarabic}
  \item The Chang--Skjelbred sequence is exact.
  \item \(\HT^{*}(X)\) is a reflexive \(R\)-module.
  \end{enumarabic}
  \goodbreak
  If \(X\) satisfies Poincaré duality, then they are also equivalent to:
  \begin{enumarabic}[resume]
  \item The equivariant Poincaré pairing
    \begin{equation*}
      \HT^{*}(X) \times \HT^{*}(X) \to R,
      \quad
      (\alpha,\beta) \mapsto \bigpair{\alpha\cup\beta,[X]} 
    \end{equation*}
    is perfect.
  \end{enumarabic}
\end{theorem}

We give two examples of \(T\)-manifolds whose equivariant cohomology is reflexive, but not free.

\begin{example}[{\cite[Sec.~6.1]{AlldayFranzPuppe:orbits1}}]
  Let \(T=(S^{1})^{n}\), and let \(N\) and~\(S\) be the north and south pole of the sphere~\(S^{2}\), respectively.
  Then
  \begin{equation}
    X = (S^{2})^{n} \setminus \{(N,\dots,N),(S,\dots,S)\}
  \end{equation}
  is a non-compact orientable \(T\)-manifold (in fact, a smooth toric variety).\footnote{%
    Strictly speaking, \(X\) does not satisfy our assumptions on \(T\)-spaces. However,
    if we instead removed small open balls around the two fixed points, we would get a \(T\)-manifold with boundary that
   satisfies our assumptions and is equivariantly homotopy-equivalent to~\(X\).}
  Its equivariant cohomology is reflexive for~\(n\ge3\), but never free.
  The latter claim can be verified by comparing Betti sums as in~\eqref{eq:betti-sum}:
  We have
  \begin{equation*}
    \dim H^{*}(X_{0}) = 2^{n}-2 < 2^{n}-1 = \dim H^{*}(X).
  \end{equation*}  
\end{example}

Compact orientable \(T\)-manifolds with reflexive, but not free equivariant cohomology
are much harder to find. According to~\cite[Cor.~1.4]{AlldayFranzPuppe:orbits1}, this is only possible for~\(n\ge5\).
All known examples are variations of the following construction,
see~\cite{Franz:2015} and~\cite{FranzHuang:2020}.

\begin{example}
  \def\uu{\mathbfit{u}}
  \def\zz{\mathbfit{z}}
  The \(S^{1}\)-action on~\(S^{3}\) given by~\(g(z,u)=(gz,u)\) for~\(g\in S^{1}\)
  and \((z,u)\in S^{3}\subset\C^{2}\) induces an action of \(T=(S^{1})^{n}\) on \((S^{3})^{n}\).
  The \newterm{big polygon space}
  \begin{equation}
    X = \bigl\{\, ((z_{1},u_{1}),\dots,(z_{n},u_{n}))\in (S^{3})^{n} \bigm| u_{1}+\dots+ u_{n} = 0 \,\bigr\}
  \end{equation}
  is a compact orientable \(T\)-manifold provided that \(n=2m+1\) is odd.
  Its equivariant cohomology is reflexive for~\(m\ge2\), but not free over~\(R\).

  The fixed point set is
  \begin{equation*}
    X_{0} = \bigl\{\, (u_{1},\dots,u_{n})\in(S^{1})^{n} \bigm| u_{1}+\dots+ u_{n} = 0 \,\bigr\}.
  \end{equation*}
  The diagonal circle action on it is free. The quotient \(X_{0}/S^{1}\) is called
  an equilateral \newterm{planar polygon space}. See~\cite{FarberSchuetz:2007} and the references given therein
  for more about these spaces.
\end{example}

\begin{remark}
  \Cref{thm:reflexive} is a special case of a more general result due to Allday--Franz--Puppe
  about the augmented \newterm{Atiyah--Bredon sequence}
  \begin{multline}
    \label{eq:ab-seq}
    \let\to\longrightarrow
    0 \to \HT^{*}(X) \stackrel{\iota^{*}}{\to} \HT^{*}(X_{0}) \to \HT^{*+1}(X_{1},X_{0}) \to \\
    \let\to\longrightarrow
    \HT^{*+2}(X_{2},X_{1}) \to \dots \to \HT^{*+n}(X_{n},X_{n-1}) \to 0,
  \end{multline}
  where \(X_{i}=\{x\in X|\rank T_{x}\ge n-i\}\) is the equivariant \(i\)-skeleton,
  and the maps are the connecting homomorphisms in the long exact sequences
  for the triples~\((X_{i+1},X_{i},\allowbreak X_{i-1})\). (In other words, apart from the leading term~\(\HT^{*}(X)\),
  this is the first page of the spectral sequence associated to the filtration~\((X_{i})\) and converging to~\(\HT^{*}(X)\).)

  The more general result alluded to characterizes the exactness of a front piece
  of the sequence~\eqref{eq:ab-seq} in terms of the syzygy order of the \(R\)-module~\(\HT^{*}(X)\).
  Syzygies are a notion from commutative algebra that allows to interpolate between torsion-free and free modules.
  (See~\cite[Sec.~1]{AlldayFranzPuppe:orbits1} for a precise definition.)
  The torsion-free \(R\)-modules are the first syzygies, the reflexive ones the second syzygies and the
  free ones the \(n\)-th syzygies.
  According to Allday--Franz--Puppe~\cite[Thm.~1.1]{AlldayFranzPuppe:orbits1},
  the sequence~\eqref{eq:ab-seq} is exact at the first~\(j\) positions (not counting the leading~\(0\))
  if and only if \(\HT^{*}(X)\) is a \(j\)-th syzygy.

  The case~\(j=1\) says that \(\HT^{*}(X)\) injects into~\(\HT^{*}(X_{0})\) if and only if it is torsion-free,
  which is a well-known consequence of the localization theorem in equivariant cohomology. The case~\(j=2\)
  is \Cref{thm:reflexive}. The case~\(j=n\) (together with a small additional argument) means that the whole sequence~\eqref{eq:ab-seq}
  is exact if and only if \(\HT^{*}(X)\) is free over~\(R\). A statement analogous to the ``if'' direction
  was proven by Atiyah in equivariant \(K\)-theory
  (apparently even preceding the work of Chang--Skjelbred) and translated to equivariant cohomology by Bredon,
  see \Cref{sec:other}.
\end{remark}

\section{Non-abelian groups}

We continue to use real coefficients.

Let \(G\) be a compact connected Lie group of rank~\(n\). Then \(R=H^{*}(BG)\) is a polynomial ring in \(n\)~indeterminates
of positive even degrees. Let \(X\) be a ``sufficiently nice'' \(G\)-space, for example a compact smooth \(G\)-manifold,
possibly with boundary, or a complex  algebraic variety with an action of the complexification of~\(G\).
In order to get a Chang--Skjelbred lemma for~\(X\), one has to use the
definitions of~\(X_{0}\) and~\(X_{1}\) that involve the ranks of the isotropy groups. We therefore set
\begin{align}
  X_{0} &= \{\, x\in X \bigm| \rank G_{x} = n \,\}, \\
  X_{1} &= \{\, x\in X \bigm| \rank G_{x} \ge n-1 \,\}.
\end{align}
Note that \(X_{0}\) contains the fixed point set~\(X^{G}\), but is larger in general.
Moreover, like~\(X_{1}\), it usually has singularities even for smooth actions.

We consider the analogue of the Chang--Skjelbred sequence in this context,
\begin{equation}
  \label{eq:cs-nonabelian}
  0 \longrightarrow \HG^{*}(X) \stackrel{\iota^{*}}{\longrightarrow} \HG^{*}(X_{0}) \stackrel{\delta}{\longrightarrow} \HG^{*+1}(X_{1},X_{0}).
\end{equation}

\begin{theorem}[Goertsches--Mare~{\cite[Thm.~2.2]{GoertschesMare:2014}}]
  Let \(X\) be a compact smooth \(G\)-manifold. If \(\HG^{*}(X)\) is free over~\(R\), then
  the Chang--Skjelbred sequence~\eqref{eq:cs-nonabelian} is exact.
\end{theorem}

\begin{theorem}[{\cite[Thm.~1.4]{Franz:2016}}]
  The Chang--Skjelbred sequence~\eqref{eq:cs-nonabelian} is exact if and only if \(\HG^{*}(X)\)
  is a reflexive \(R\)-module.
\end{theorem}

Using induction from a maximal torus~\(T\subset G\) to~\(G\),
one can adapt the examples of spaces with reflexive torus-equivariant cohomology
given in the previous section to actions of non-abelian groups,
see~\cite[Prop.~4.5, Example~6.13]{Franz:2016}.

\section{Other generalizations}
\label{sec:other}

We briefly mention other generalizations of the Chang--Skjelbred lemma.

Allday--Franz--Puppe~\cite{AlldayFranzPuppe:orbits2} consider actions of \(p\)-tori~\((\Z/p)^{n}\).
The results are largely analogous to the torus case, with some differences between the cases~\(p=2\) and \(p\)~odd.
The proofs given in~\cite{AlldayFranzPuppe:orbits2} are new and work also for tori and other compact connected Lie groups.
An alternative proof for the case~\(p=2\) has appeared in recent work of Bourguiba--Lannes--Schwartz--Zarati~\cite{BLSZ}.
Their approach is based on the action of the Steenrod algebra.

Goertsches--Töben~\cite[Thm.~A]{GoertschesToeben:2010} have generalized the Chang--Skjelbred lemma to torus actions without fixed points.
The freeness of~\(\HT^{*}(X)\) is replaced by the condition that the equivariant cohomology be Cohen--Macaulay.
The role of the fixed point set~\(X_{0}\) is taken by~\(X_{m}\) where \(m\) is the smallest orbit dimension occurring in~\(X\).

We finally discuss the analogue of the Chang--Skjelbred lemma for torus-equiv\-ari\-ant \(K\)-theory
(with compact supports),
\begin{equation}
  \label{eq:cs-ktheory}
  0 \longrightarrow \KT^{*}(X) \stackrel{\iota^{*}}{\longrightarrow} \KT^{*}(X_{0}) \stackrel{\delta}{\longrightarrow} \KT^{*+1}(X_{1},X_{0}).
\end{equation}
Atiyah~\cite{Atiyah:1974} proved the exactness of this sequence under the assumption that \(\KT^{*}(X)\)
is free over the representation ring~\(R=\Z[T]=\Z[t_{1},t_{1}^{-1},\dots,t_{n},t_{n}^{-1}]\).
In fact, he proved the exactness of the \(K\)-theoretic analogue
\begin{multline}
  \label{eq:ab-seq-kt}
  \let\to\longrightarrow
  0 \to \KT^{*}(X) \stackrel{\iota^{*}}{\to} \KT^{*}(X_{0}) \to \KT^{*+1}(X_{1},X_{0}) \to \\
  \let\to\longrightarrow
  \KT^{*+2}(X_{2},X_{1}) \to \dots \to \KT^{*+n}(X_{n},X_{n-1}) \to 0
\end{multline}
of the Atiyah--Bredon sequence~\eqref{eq:ab-seq}.
After the publication of Goresky--Kottwitz--MacPherson's paper~\cite{GoreskyKottwitzMacPherson:1998},
Rosu--Knutson~\cite{KnutsonRosu:2003} gave a proof of the exactness of the sequence~\eqref{eq:cs-ktheory}, tensored with~\(\C\).

\begin{remark}
  We elaborate on how Atiyah's results give the exact sequence~\eqref{eq:ab-seq-kt}.
  In~\cite[Lecture~7]{Atiyah:1974}, Atiyah works with a complex \(T\)-module~\(E\), but his arguments
  work equally well for the \(T\)-spaces~\(X\) we consider. Assuming that \(\KT^{*}(X)\) is free over~\(R\),
  he proves that for any~\(0\le i\le n\) there is a short exact sequence
  \begin{equation}
    \let\to\longrightarrow
    0 \to \KT^{*}(X,X_{i-1}) \to \KT^{*}(X_{i},X_{i-1}) \stackrel{\delta}{\to} \KT^{*+1}(X,X_{i}) \to 0.
  \end{equation}
  see~\cite[eq.~(7.3)]{Atiyah:1974}. (Here \(X_{-1}=\emptyset\).  Also note that Atiyah uses
  the notation~\(M(i)=\KT^{*}(X,X_{n-i})\) and~\(L(i)=\KT^{*}(X_{n-i},X_{n-i-1})\).)
  As explained in~\cite[Lemma~4.1]{FranzPuppe:2007},
  these short exact sequences can be spliced together to give the exact sequence~\eqref{eq:ab-seq-kt}.
  The latter appeared first in a paper by Bredon~\cite[Main Lemma]{Bredon:1974}, in the translation to \(T\)-equivariant cohomology.
\end{remark}

\end{document}